\documentclass[10pt]{amsart}
\usepackage{amsmath, amssymb, amscd, epsf}

\newcommand{\D}{\ensuremath{\mathcal{D}}}
\newcommand{\PP}{\ensuremath{\mathbb{P}}}
\newcommand{\ZZ}{\ensuremath{\mathbb{Z}}}
\newcommand{\NN}{\ensuremath{\mathbb{N}}}
\newcommand{\FF}{\ensuremath{\mathbb{F}}}

\newcommand{\RR}{\ensuremath{\mathbb{R}}}

\newcommand{\ba}{{\bf a}}
\DeclareMathOperator{\reg}{reg}
\DeclareMathOperator{\rank}{rank}
 
\DeclareMathOperator{\pdim}{proj-dim}
\DeclareMathOperator{\sign}{sign} \DeclareMathOperator{\lcm}{lcm}
\DeclareMathOperator{\tor}{Tor} \DeclareMathOperator{\pos}{pos}

\DeclareMathOperator{\udeg}{deg_{\bf u}}
\DeclareMathOperator{\vdeg}{deg_{\bf v}}
\DeclareMathOperator{\vreg}{reg_{\bf v}}
\DeclareMathOperator{\regnum}{reg-num}
\DeclareMathOperator{\vregnum}{reg-num_{\bf v}}
\DeclareMathOperator{\uregnum}{reg-num_{\bf u}}

\newtheorem{thm}{Theorem}[section]
\newtheorem{cor}[thm]{Corollary}
\newtheorem{lem}[thm]{Lemma}

\theoremstyle{definition}

 \newtheorem{defin}[thm]{Definition}
\newtheorem{ex}[thm]{Example}

\newtheorem{rmk}[thm]{Remark}

\begin{document}

\title{Multigraded regularity: coarsenings and resolutions}
\author{Jessica Sidman}
\address{Department of Mathematics and Statistics \\
451A Clapp Lab\\
Mount Holyoke College\\
South Hadley, MA 01075, USA}
\email{jsidman@mtholyoke.edu}
\author{Adam Van Tuyl}
\address{Department of Mathematical Sciences \\
Lakehead University \\
Thunder Bay, ON P7B 5E1, Canada}
\email{avantuyl@sleet.lakeheadu.ca}
\author{Haohao Wang}
\address{Department of Mathematics\\
Southeast Missouri State University \\
Cape Girardeau, MO 63701, USA}
\email{hwang@semo.edu}
\thanks{Last updated: September 19, 2005}
\begin{abstract}
Let $S = k[x_1,\ldots,x_n]$ be a $\ZZ^r$-graded ring  with
$\deg (x_i) = {\bf a}_i \in \ZZ^r$ for each $i$ and suppose that $M$ is  a finitely generated
$\ZZ^r$-graded $S$-module.  In this paper we describe how to find finite subsets
of $\ZZ^r$ containing the multidegrees of the minimal multigraded syzygies of $M.$  To find such a set, we first
coarsen the grading of $M$ so that we
can view $M$ as a $\ZZ$-graded $S$-module.
We use a generalized notion
of Castelnuovo-Mumford regularity, which was introduced  by D. Maclagan and G. Smith,
to associate to $M$ a number which we call the regularity number of $M$.  The
minimal degrees of the multigraded minimal syzygies are bounded in terms
of this invariant.
\end{abstract}

\maketitle


\section{Introduction}
Let $S = k[x_1,\ldots, x_n]$ with $\deg (x_i) = 1$ and $M$ be a finitely generated graded $S$-module.  The Castelnuovo-Mumford regularity of $M,$
denoted $\reg(M),$ is a cohomological invariant that bounds the ``size'' of
its minimal free resolution.
The module $M$ is $d$-regular if
$H^i_{\mathfrak{m}}(M)_p = 0$ for all $i \geq 0$ and all $p \geq d-i+1,$ where
$\mathfrak{m} = \langle x_1, \ldots, x_n \rangle.$  The \emph{regularity} of $M$ is the smallest integer $d$
for which $M$ is $d$-regular.
 If
\[ 0 \to F_r \to \cdots \to F_i \to \cdots \to F_0 \to M \to 0\]
is a minimal free graded resolution of $M,$ then the degrees of the generators
of $F_i$ are bounded above by $ \reg(M)+i.$
Hence, $\reg(M)$ gives
a finite set which contains all of the possible degrees of minimal generators of the modules $F_i.$

The goal of this paper is to find a cohomological theory of regularity
which will complement that of \cite{MS} and give finite bounds on the degrees of the minimal syzygies of $M$ in a $\ZZ^r$-graded setting.
More precisely, suppose that $S = k[x_1,\ldots,x_n]$ is $\ZZ^r$-graded with
$\deg (x_i) = {\bf a}_i \in \ZZ^r$ for each $i$, and let $M$ be a finitely generated
$\ZZ^r$-graded $S$-module.
Under suitable hypotheses on the grading, $M$ has a finite $\ZZ^r$-graded minimal free resolution
(cf. Definition \ref{defn: pm}).  We say that a finite subset $D$
of $\ZZ^r$ such that the multidegree of every minimal generator of $M$ belongs to $D$ is a \emph{finite bound} on the degrees of the generators of $M.$  

Partial solutions to this problem have appeared in previous investigations into the regularity
of multigraded modules (cf. \cite{ACD,H,HW,MS,SV}).  We summarize these approaches here.

Maclagan and Smith \cite{MS} developed a multigraded theory of regularity for
sheaves on a simplicial toric variety $X$ with an algebraic variant
defined in terms of the vanishing of graded pieces of
$H_B^i(M),$ the $i$th local
cohomology module of $M$. Here $B$ is the irrelevant ideal of the
homogeneous coordinate ring of $X$
which is multigraded by a finitely generated
abelian group $A.$
Multigraded regularity, which we denote by $\reg_B(M)$ and refer to as the
\emph{B-regularity} of $M,$
is a subset of $A,$ and the degrees of the minimal syzygies of $M$ lie in
the complement of certain shifts of $\reg_B(M).$  As we see in Example
\ref{ex: points}, these complements may not
be bounded if $\rank(A)$ $\geq 2,$ so $B$-regularity does not necessarily give
 a finite bound on the degrees of the minimal syzygies.

\begin{ex}\label{ex: points}
Consider the bigraded coordinate ring of $\PP^1 \times \PP^1$, that is, $S = k[x_0,x_1,y_0,y_1]$
with $\deg (x_i) = (1,0)$ and $\deg (y_i) = (0,1)$. Let $I$ be the defining ideal of
the set of points
\[X = \{[1:0] \times [1:0],[1:0] \times [0:1],[0:1] \times [1:0], [0:1] \times [0:1]\}.\]
The bigraded Hilbert function of $S/I$ is
\[
H_{S/I} =
\left[\begin{matrix}
\vdots&\vdots&\vdots&\\
2 & 4 & 4 & \cdots\\
2 & 4 & 4 & \cdots\\
1 & 2 & 2 & \cdots
\end{matrix}\right]
\]
where we identify the value at the $(i,j)$th integer cartesian coordinate with the value of $H_{S/I}(i,j)$.
Applying Proposition 6.7 of \cite{MS} gives
\[\reg_B(S/I) = \{(i,j) \in \NN^2 ~|~ H_{S/I}(i,j) = \deg X = 4\} = (1,1) + \NN^2.\]
The degrees of the minimal generators are contained in the unbounded set
\[\NN^2 \backslash \left(((1,2) + \NN^2)  \cup ((2,1) + \NN^2) \right).\]
More generally, the degrees of the minimal $i$th syzygies are
contained in
\[\NN^2 \backslash \bigcup_{m,n \geq 0, m+n = i} ((1,1) + (m,n)) + \NN^2,\]
which is not finite.  To foreshadow the approach we shall use
in this paper, note that the ideal $I = \langle x_0x_1,y_0y_1 \rangle$ is a $\ZZ$-graded complete intersection with $\reg(S/I) = 2$.  So,
if $F$ is any minimal $i$th syzygy with $\deg F = (a,b)$, then $a + b \leq 2 + i$.
There are only a finite number of $(a,b) \in \NN^2$ with this property.
\end{ex}

Hoffman and the third author \cite{HW} introduced a notion of
strong regularity for bigraded modules over the homogeneous coordinate ring of
$\PP^{n} \times \PP^{m},$
$S=k[x_0,\ldots,x_n,y_0,\ldots,y_m],$  with $\deg (x_i) = (1,0)$
and $\deg (y_i) = (0,1).$   Strong regularity, which is defined in
terms of the vanishing of graded pieces of the local cohomology modules $H_{{\bf
x}}^i(M), H_{\bf y}^i(M),$ and   $H_{\bf x + y}^i(M)$ where ${\bf
x}=\langle x_1,\ldots,x_n\rangle$ and ${\bf y} = \langle
y_1,\ldots,y_m\rangle$, gives a finite set that bounds the degrees
of the minimal syzygies of $M.$  This approach relies on the
Mayer-Vietoris sequence to relate the various local cohomology
modules. Generalizations to this approach may be possible using a
more complicated spectral sequence, but we take a simpler approach
here.

Extending earlier work of Aramova, Crona, and De Negri \cite{ACD} (see also R\"omer \cite{R}), the first
two authors \cite{SV}
studied the resolutions of finitely generated modules over the coordinate ring
of $\PP^{n_1} \times \cdots \times \PP^{n_r}$.  If $M$ is such a module, one can assign to $M$ a vector
$\underline{r}(M) =(d_1,\ldots,d_r) \in \NN^r$ called the
resolution regularity vector of $M$.  This vector provides a finite bound on
the degrees of the minimal syzygies of $M$; more precisely, if $F$ is a minimal multigraded $i$th syzygy
of $M$, then $\underline{0} \leq \deg F \leq (d_1+i,\ldots,d_r+i)$.
Work of H\`a in \cite{H} and results in \cite{SV} investigate the relationship between
the resolution regularity vector and the definition of Maclagan and Smith.  H\`a also introduced
a multigraded analog of the $a^*$-invariant to study both $\underline{r}(M)$ and $\reg_B(M)$.

Our point of view is shaped by the observation that in the work of
Hoffman and the third author a finite bound on the degrees of syzygies
may be obtained by
considering only vanishings of $H_{\bf x + y}^i(M),$ and this may be
extended to nonstandard multigraded rings.
We also exploit the fact that Maclagan and Smith's
definition of regularity does provide a finite bound on the degrees of syzygies if
$M$ is a finitely generated module over the coordinate ring of a weighted
projective space, that is, the nonstandard $\ZZ$-graded ring $S = k[x_1,\ldots,x_n]$
with $\deg(x_i) = a_i \in \ZZ$ for each $i$.

Our method can be summarized as follows:  Let $M= \bigoplus_{{\bf a} \in \ZZ^r} M_{\bf a}$
be a finitely generated $\ZZ^r$-graded module over a positively multigraded (see
Definition \ref{defn: pm}) $\ZZ^r$-graded ring $S=k[x_1,\ldots,x_n].$
We can coarsen the grading of $M$ to form a $\ZZ$-graded module
by picking a suitable vector ${\bf v}$ and setting $M^{[{\bf v}]} =
\bigoplus_{m\in \ZZ} (\bigoplus_{{\bf a}\cdot{\bf v}=m} M_{\bf a})$.
The module $M^{[{\bf v}]}$ is an $S^{[{\bf v}]}$-module where $S^{[{\bf v}]}$ is the
nonstandard graded polynomial ring with $\deg(x_i) = {\bf a}_i \cdot {\bf v}$.
We associate to $M$ a {\it regularity number}, denoted $\vregnum(M)$, where
$\vregnum(M)$ is defined using the regularity of $M^{[{\bf v}]}$ as defined by Maclagan and Smith.  Our
main result (cf. Corollary \ref{cor: minimal_generator}) uses $\vregnum(M)$ to
obtain finite bounds on which multidegrees may appear in a minimal $\ZZ^r$-graded
free resolution of $M$.

Different vectors {\bf w} result in different coarsenings, and
hence, $\regnum_{\bf w}(M)$  may give us different bounds on the
multidegrees. We examine how we can use different coarsenings of
$M$ to improve our bounds on the multidegrees. Even though there
are an infinite number of possible coarsening vectors {\bf w} (and
hence, different degree bounds), we show (cf. Theorem \ref{thm:
mspcv}) that all of the information on degree bounds that can be
obtained from coarsenings can be obtained from a finite number of
such vectors. We call a finite set of vectors with this property a
minimal set of
positive coarsening vectors.

This paper is structured as follows.  In \S 2 we recall important
notions related to multigraded rings and positive coarsenings of
the grading.  In \S 3  we introduce the idea of the regularity
number of a $\ZZ^r$-graded module $M$ under a given positive
coarsening, and we show how this number can be used to bound
degrees appearing in a minimal free resolution of $M.$
In \S 4 we show that scalar multiples of a positive coarsening vector all give the same degree bounds.
  In \S 5 we relate $\vregnum(M)$ to the notion of the multigraded regularity in \cite{MS} and the notion of
regularity based on resolutions in \cite{SV}.  In \S 6 we illustrate the theory
that has come before with examples.

\noindent{\bf Acknowledgments} We thank David Cox for suggesting
the collaboration and for his encouragement, and T\`ai H\`a, Robert Lazarsfeld, Diane Maclagan and Greg Smith for their comments
and suggestions. We computed many
examples using the computer algebra packages \emph{CoCoA}
\cite{Co} and \emph{Macaulay 2} \cite{GS}.  The first author
was partially supported by an NSF postdoctoral fellowship at the University
of Massachusetts in Amherst in 2004-2005 on Grant No. 0201607 and also thanks the Clare Boothe Luce Foundation for support.
The second author thanks NSERC for financial support while working
on this project and the hospitality of Mt. Holyoke College while visiting
the first author.  The third author thanks NSF-AWM and GRFC for
travel support.


\section{Preliminaries}

We briefly summarize the relevant
definitions and properties of multigraded rings and multigraded
resolutions in \S 2.1.  We refer the reader to Chapter 8 of \cite{Mi-S}
for a comprehensive introduction to multigraded
rings and modules.  We describe the notion of coarsening
a multigrading in \S 2.2.  (See also Chapters 7 and 8 of \cite{Mi-S}.)

\subsection{Multigraded polynomial rings and resolutions}
\subsubsection{Multigraded rings and modules}
Let $S = k[x_1,\ldots,x_n]$ with a $\ZZ^r$-grading given by
 a group homomorphism $\deg: \ZZ^n
\rightarrow \ZZ^r$ where
$\deg(x_i) = \ba_i$ for some $\ba_i \in \ZZ^r$.  Let $Q = \deg(\NN^n)$ denote the
subsemigroup of $\ZZ^r$ generated by ${\bf a}_1,\ldots,{\bf a}_n$.

If $S$ has a $\ZZ^r$-grading, then we say $F \in S$ is homogeneous of degree $\ba$
if all the terms of $F$ have degree $\ba$.
We let $S_{\ba}$ denote the $k$-vector space consisting of all the homogeneous
forms of degree $\ba$ for each $\ba \in \ZZ^r$.
A finitely generated $S$-module $M$ is $\ZZ^r$-graded if
\[M = \bigoplus_{ {\bf a} \in \ZZ^r} M_{\bf a} ~~\mbox{and}~~
S_{\bf a}M_{\bf b} \subseteq M_{{\bf a}+{\bf b}} ~~\mbox{for all ${\bf a},{\bf b} \in \ZZ^r$}.\]
For any ${\bf a} \in \ZZ^r,$
define $M({\bf a})$ to be the finitely generated $\ZZ^r$-graded $S$-module where
$M({\bf a})_{\bf p} = M_{\bf p+a}$ for all ${\bf p}$.

\subsubsection{Multigraded resolutions}
Given a finitely generated $\ZZ^r$-graded $S$-module $M$, we shall be interested in the
$\ZZ^r$-graded minimal free resolution of $M$.  However, this notion may not be well-defined because the notion of minimality breaks down without additional hypotheses.
To obtain the desired behavior, we will impose the additional constraints defined below
on the grading of $S$ that guarantee a version of Nakayama's Lemma.  (See Chapter 8.2 in \cite{Mi-S}.)

\begin{defin}[Definition 8.7 in \cite{Mi-S}]\label{defn: pm}
The polynomial ring $S$ is \emph{positively multigraded} by $\ZZ^r$ if $\deg(x_i)
\neq 0$ for all $i$ and the semigroup $Q$ defined above has no nonzero invertible
elements.
\end{defin}

The Cox homogeneous coordinate ring of a complete smooth toric variety is an important example of a positively multigraded ring.  We will say a bit more about this case in \S \ref{section: B-reg}.  Example \ref{ex: proj prod} is a special case of such a ring.

\begin{ex}\label{ex: proj prod}
We say that the polynomial ring $S$ is a {\it standard multigraded ring} if
$S$ is the homogeneous coordinate ring of the product of projective spaces
$\mathbb{P}^{n_1} \times \cdots \times \mathbb{P}^{n_r}$.
Equivalently, $S = k[x_{1,0}, \ldots, x_{1,n_1}, \ldots, x_{r,0}, \ldots, x_{r,n_r}]$ is
$\ZZ^r$-graded by setting $\deg (x_{i,j}) = e_i$, the $i$th standard basis vector.
\end{ex}

We assume throughout this paper that
$S$ is positively multigraded.
Consequently, if $M$ is a finitely generated
$\ZZ^r$-graded $S$-module, then $\dim_k M_{\bf a} < \infty$ for all ${\bf a} \in \ZZ^r.$
(See Theorem 8.6 in \cite{Mi-S} for a proof and other equivalent properties.)

In this  setting there is a well-defined notion of a minimal free resolution of
a $\ZZ^r$-graded module $M$ (cf. Chapter 8.3 \cite{Mi-S}).  That is, there
exists an exact complex of the form
\begin{equation}\label{multigraded resolution}
{\bf F}_{\bullet} :0 \rightarrow F_{\ell} \rightarrow F_{\ell-1}
\rightarrow \cdots \rightarrow F_1 \rightarrow F_0 \rightarrow M
\rightarrow 0
\end{equation}
such that $\ell \leq n$, and $F_i = \bigoplus_{\ba \in \ZZ^r} S(-\ba)^{\beta_{i,\ba}(M)}$ is finitely generated.
The number
\[\beta_{i,\ba}(M)= \dim_k \tor_i^S(M,k)_{\bf a}\] is the $i$th graded
 Betti number of $M$ of degree ${\bf a}$. This invariant of $M$ counts the number
of minimal generators of degree $\ba$ in
the $i$th syzygy module of $M$.
Our goal is to find a mechanism for finding a finite set $\mathcal{D}_i(M) \subset \ZZ^r$
such that  $\beta_{i,\ba}(M)$ is always zero for $\ba \in \ZZ^r \backslash \mathcal{D}_i(M)$.  In
other words, we wish to provide a finite set of possible values for the degrees
of the minimal $i$th syzygies.

Although the degrees of elements of $S$ lie in $Q$,
the degrees of elements of $M$ may not.  However,
 by Theorem 8.20 of \cite{Mi-S}, the
multigraded Hilbert series of $M,$
$H(M; {\bf t}) = \sum_{\ba \in \ZZ^r} \dim_k (M_{\bf a}){\bf t^a}$
lies in
$\ZZ[[Q]][\ZZ^r] = \ZZ[[Q]] \otimes_{\ZZ[Q]} \ZZ[\ZZ^r],$ the ring of Laurent series supported on \emph{finitely
many translates of $Q$} (see pg.155 of
\cite{Mi-S}). This
just means that the degrees of elements of $M$ lie in
$({\bf b}_1 + Q) \cup \cdots \cup  ({\bf b}_k +Q)$
for some finite collection of ${\bf b}_i$s in $\ZZ^r.$

\subsection{Coarsening Gradings}
If $M$ is a finitely generated $\ZZ^r$-graded $S$-module, we  can
\emph{coarsen} the $\ZZ^r$-grading, thus allowing us to view $M$
as a $\ZZ$-graded module (see Corollary 7.23 in \cite{Mi-S}).
We pass to a $\ZZ$-grading by choosing a vector ${\bf v} \in \ZZ^r$ and defining
$\vdeg (m) := \deg (m) \cdot {\bf v}$ for $m \in M$ using the dot product
in $\ZZ^r.$  Note that
$ \vdeg(x_i m) = \vdeg(x_i) + \vdeg(m)$
for all variables $x_i$ and elements $m \in M.$

We write $S^{[{\bf v}]}$ to denote the $\ZZ$-graded
polynomial ring $S^{[{\bf
v}]}$ with $(S^{[{\bf v}]})_m = \bigoplus_{\ba\cdot{\bf v} = m}
S_{\ba}$ and $M^{[{\bf v}]}$ for the
$\ZZ$-graded $S^{[{\bf v}]}$-module
$M^{[{\bf v}]} = \bigoplus_{m \in \ZZ} \left(\bigoplus_{{\bf a}\cdot{\bf v} =m} M_{\bf a} \right).$
If the coarsening vector is clear, we will sometimes drop the superscript.

We want the $\ZZ$-graded ring $S^{{[\bf v}]},$ which can be viewed
as the coordinate ring of a weighted projective space, to be
positively graded.

\begin{defin}
A vector ${\bf v} \in \ZZ^r$ is a \emph{positive coarsening vector} for the
$\ZZ^r$-graded ring $S$ if $\vdeg(x_i) > 0$ for all $i.$
\end{defin}

\begin{lem}
If $S$ is positively $\ZZ^r$-graded, then a positive coarsening vector exists.
\end{lem}

\begin{proof}
The set
\[ \pos(Q) = \left.\left\{ \sum \lambda_i {\bf a}_i ~\right|~ \lambda_i \in \RR_{\ge 0} \right\}\]
is the convex cone in $\RR^r$ generated by the vectors in $Q.$  Our assumption
that $S$ is positively multigraded
implies that $\pos(Q)$ is a pointed cone, i.e., $\pos(Q)$ does not
contain any nontrivial linear subspaces.
If $\pos(Q)$ is pointed, then ${\bf 0}$ is a face, so there
exist hyperplanes $H \subset \RR^r$ so that $\pos(Q) - \{ {\bf 0} \}$ lies
strictly on one side of $H$.  If we take ${\bf v}$ to be a normal vector to $H$
on the same side of $H$ as $\pos(Q),$ then $\vdeg(x_i) > 0$ for each $i.$
\end{proof}

\begin{rmk}
Let $S$ be positively multigraded by $\ZZ^r$ with degrees in
the subsemigroup $ Q\subset \ZZ^r$, and let $M$ be a finitely generated
$\ZZ^r$-graded $S$-module.  If ${\bf v} \in \ZZ^r$ is a positive coarsening
vector, then $\dim_k (M^{[{\bf v}]})_m <\infty$ for
all $m \in \ZZ$ since the ring $S^{[{\bf v}]}$ is positively graded.
\end{rmk}

Suppose ${\bf v}$ is a positive coarsening vector for $S$, and that
we have a multigraded minimal free resolution of
the $\ZZ^r$-graded $S$-module $M$ of the form $(\ref{multigraded resolution})$.  Using
${\bf v}$ we can pass to a $\ZZ$-graded minimal free
resolution of $M^{[{\bf v}]}$ as an $S^{[{\bf v}]}$-module:
\begin{equation}
0 \rightarrow F^{[{\bf v}]}_{\ell} \rightarrow
 F^{[{\bf v}]}_{\ell-1} \rightarrow \cdots
\rightarrow F^{[{\bf v}]}_{1} \rightarrow
 F^{[{\bf v}]}_{0} \rightarrow M^{[{\bf v}]} \rightarrow 0
\end{equation}
where
\begin{eqnarray*}
F^{[{\bf v}]}_i &= &\bigoplus_{\ba \in \ZZ^r} S^{[{\bf v}]}(-\ba\cdot {\bf v})^{\beta_{i,{\ba}}(M)}
 =  \bigoplus_{m \in \ZZ} S^{[{\bf v}]}(-m)^{\sum_{\ba\cdot{\bf v} = m} \beta_{i,{\ba}}(M)}.
\end{eqnarray*}
The above expression allows us to deduce that
$\sum_{\ba\cdot{\bf v} = m} \beta_{i,{\ba}}(M) = \beta_{i,m}(M^{[{\bf v}]}).$
The heart of our method is to find bounds on $m$ such that  $\beta_{i,m}(M^{[{\bf v}]}) \neq 0$
by investigating the regularity of $M^{[{\bf v}]}$ as a $\ZZ$-graded $S^{[{\bf v}]}$-module.  This, in turn,
gives us information on which ${\bf a} \in \ZZ^r$ have the property that $\beta_{i,{\bf a}}(M) \neq 0$.


\section{Regularity of $\ZZ$-graded modules}

Our main tool will be a special case of the multigraded regularity
introduced in \cite{MS} which is also related to the commutative version of the notion
of regularity introduced in \cite{B}.
Throughout this section we assume that $M$ is a finitely
generated $\ZZ^r$-graded $S$-module and that ${\bf v} \in \ZZ^r$ is a positive
 coarsening vector.

Specializing Definition 4.1 in \cite{MS} gives us:
\begin{defin}
Let $M$ be a finitely generated $\ZZ^r$-graded $S$-module $M$
and let ${\bf v}$ be a positive coarsening vector. Let
$\mathfrak{m} = \langle x_1, \ldots, x_n \rangle$ and $c_{\bf v}=
\lcm(\vdeg(x_i))_{i=1}^n.$  Define
\[\vreg(M) = \{ p \in \ZZ \mid H^i_{\mathfrak{m}}(M)_q = 0 ~~\forall i \geq 0, ~~\forall q \in p + \NN c_{\bf v}[1-i]\}\]
where
\[\NN c_{\bf v}[i] = ic_{\bf v} + \NN c_{\bf v} = \{c_{\bf v}(i + d) ~|~ d \in \NN\}.\]
\end{defin}
Note that $\vreg(M)$ is a special case of the definition of the regularity of $M^{[{\bf v}]}$ given in \cite{MS}.

\subsection{The regularity number}

The vanishing conditions required by the definition of $\vreg(M)$
are only required for shifted multiples of $c_{\bf v}$.  So, if $c_{\bf v} >1,$ then $p \in \vreg(M)$ may not
imply that $q \in \vreg(M)$ for all $q \ge p$ as demonstrated in the following example.

\begin{ex}\label{ex: gaps}
Suppose that $S = K[x_1, x_2]$ with $\deg(x_i) = 4$ for $i=1,2$.  Then
since $H^0_{\mathfrak{m}}(S) = H^1_{\mathfrak{m}}(S) = 0,$
\begin{eqnarray*}
\reg(S) &= &\{u \in \ZZ \mid H^2_{\mathfrak{m}}(S)_p = 0 \ \forall p \in u + 4\NN[1-2]  \} \\
&=& \{u \in \ZZ \mid H^2_{\mathfrak{m}}(S)_p = 0 \ \forall p \in u -4 + 4\NN \}.
\end{eqnarray*}
If we compute $H^2_{\mathfrak{m}}(S)$ using a \v{C}ech complex, then  $H^2_{\mathfrak{m}}(S)$
is a quotient of $S_{x_1x_2}$
by all elements of $S_{x_1x_2}$ that do not have both $x_1$ and $x_2$ in the
denominator when written in lowest terms.  Since all of the elements of $S_{x_1x_2}$
have degrees that are
multiples of 4, $H^2_{\mathfrak{m}}(S)_p = 0$ for all $p$ that are not multiples of 4.

Of course it is also true that $H^2_{\mathfrak{m}}(S)_p = 0$ for all
$p > -8$ since any monomial written in lowest terms that has degree $> -8$
cannot have both $x_1$ and $x_2$ in the denominator.
  So in this example, if $q > -4,$ then
\[H^2_{\mathfrak{m}}(S)_{q-4+4m} = 0 ~~\mbox{for all $m \in \NN$}\]
since $q-4 > -8.$  Therefore,
$q\in \reg(S),$ if $q > -4$, and as well, $-5 \in \reg(S)$.   However $-4 \not\in\reg(S)$
since $H^2_{\mathfrak{m}}(S)_{-4} \neq 0$.
\end{ex}

Despite the behavior exhibited in Example \ref{ex: gaps}, we can guarantee that there is an integer $p \in \vreg(M)$
such that $q \in \vreg(M)$ for all $q\geq p.$

\begin{thm} \label{thm: vregnum}
 Let $M$ be a finitely generated $\ZZ^r$-graded $S$-module.
Then there exists a $p \in \vreg(M)$ such that $q \ge p$ implies that $q \in \vreg(M)$.
\end{thm}

In light of the previous result we have:

\begin{defin} The {\it regularity number} of $M$  (with respect to a positive coarsening
vector ${\bf v}$), denoted $\vregnum(M)$, is
\[\vregnum(M) :=\inf\{p \in \vreg(M) ~|~ q \in \vreg(M) ~\mbox{for all}~ q \geq p \}.
\]
\end{defin}

We need the following lemma to prove Theorem \ref{thm: vregnum}.

\begin{lem}\label{L:vregnums}
There exists a $p \in \vreg(S)$ such that $q \ge p$ implies that $q \in \vreg(S)$.
In particular
\[\vregnum(S) = (n-1)c_{\bf v}+1 -\sum \vdeg(x_i).\]
\end{lem}

\begin{proof}
Let $\mathfrak{m} = \langle x_1,\ldots, x_n \rangle.$  The set
\begin{eqnarray*}
\vreg(S) &= &\{ u \in \ZZ \mid H_{\mathfrak{m}}^i(S)_w = 0
\ \forall i \ge 0 \ \mathrm{and}~ \forall w \in u + \NN c_{\bf v}[1-i] \}\\
&=&  \{ u \in \ZZ \mid H_{\mathfrak{m}}^{n}(S)_w = 0 \ \forall w \in u + \NN c_{\bf v}[1-n] \}\\
&=&  \{ u \in \ZZ \mid H_{\mathfrak{m}}^{n}(S)_w = 0 \ \forall w \in u + (1-n)c_{\bf v} + \NN c_{\bf v}\}.
\end{eqnarray*}
The module  $H^{n}_{\mathfrak{m}}(S)$ has a $\ZZ^n$-grading
and
$H^{n}_{\mathfrak{m}}(S)_{{\bf q}} \neq 0$ for ${\bf q} \in \ZZ^n$
if and only if all the coordinates of ${\bf q}$ are negative (see \cite{EMS}). Therefore, it is certainly true that
$H^{n}_{\mathfrak{m}}(S)_w = 0$ if  $w > -\sum \vdeg(x_i).$
Hence
\[ \{ u \in \ZZ \mid u> (n-1)c_{\bf v} - \sum \vdeg(x_i)\} \subseteq \vreg(S).\]
Note that the inclusion may be strict as in Example \ref{ex: gaps}.

We see that $\vregnum(S)$ cannot be equal to $(n-1)c_{\bf v} - \sum \vdeg(x_i)$
since then we would have $H^n_{\mathfrak{m}}(S)_w = 0$
where
\[ w = (n-1)c_{\bf v} - \sum \vdeg(x_i)+(1-n)c_{\bf v} = -\sum \vdeg(x_i).
\]
Since $-\sum \vdeg(x_i)$ is the total degree of a Laurent monomial in which
all of the exponents are negative, we know that $H^n_{\mathfrak{m}}(S)_w \neq 0.$
\end{proof}

\begin{cor}\label{C:pdim0}Let $S$ be as above and let ${\bf d}_i \in \ZZ^r.$  Then
\[\vregnum(S(-{\bf d}_1)\oplus \cdots \oplus S(-{\bf d}_n)) = \max\{\vregnum(S)+{\bf d}_i \cdot {\bf v}\}_{i=1}^n.\]
\end{cor}

\begin{proof}
 This follows from the fact that
 \[H^i_{\mathfrak{m}}(S(-{\bf d}_1)\oplus \cdots \oplus S(-{\bf d}_n))_w \cong H^i_{\mathfrak{m}}(S)_{w-{\bf d}_1 \cdot {\bf v}}
\oplus \cdots \oplus H^i_{\mathfrak{m}}(S)_{w-{\bf d}_n \cdot {\bf v}}.\]
\end{proof}

\begin{proof}[Proof of Theorem \ref{thm: vregnum}]
We use induction on the projective dimension of $M$.  The result is trivial
if $M = 0,$ and if $M$ has projective dimension zero, then the result holds by the previous lemma since $M$ is free
in this case.

Suppose now that the projective dimension of $M$ is greater than $0.$  Then there exists a short exact
sequence
\[ 0 \to M' \to F \to M \to 0\]
where $F$ is a free module.  Since $M'$ has
projective dimension strictly less than the projective dimension of $M,$ the
result holds for $M'.$  Therefore, there exist $p_0, p_1$ such that
$p_0 \in \vreg(M')$ implies $q \in \vreg(M')$ for all $q \ge p_0$ and $p_1 \in
\vreg (F)$ implies $q \in \vreg(F)$ for all $q \ge p_1.$
Using Lemma 7.1 of \cite{MS}, we see that
\[ (\vreg(M')-c_{\bf v} )\cap \vreg(F) \subseteq \vreg(M).\]
Hence $q \in \vreg(M)$ for all $q \ge \max\{p_0-c_{\bf v}, p_1\}.$
\end{proof}

\subsection{Bounds on the degrees of the syzygies}
The following theorem and its corollary are the main results of this paper.
Both results give bounds on
the degrees in $\ZZ^r$ that may appear in a free resolution of $M.$  These results
are similar to those in \S 5 of \cite{B} and \S 7 of \cite{MS} (see
\S \ref{sec: B-reg}).

\begin{thm}\label{thm: free res}
Let $M$ be a finitely generated $\ZZ^r$-graded $S$-module, and \[s_{\bf v} := \max\{nc_{\bf v} - \sum \vdeg(x_i), c_{\bf v}\}.\]
If $F$ is a minimal generator of
the $i$th syzygy module of $M^{[{\bf v}]}$, then
\[\deg F \leq \vregnum(M)+is_{\bf v} +
c_{\bf v}-1.\]
\end{thm}

\begin{rmk}
Note that if $\vdeg(x_i)= 1$ for all $i$, then $c_{\bf v} = s_{\bf v} =1$ and the statement of the
above theorem  specializes to the statement that
$F_i$ is generated by elements of degree at most $\vregnum(M)+is_{\bf v} +1-1 = \vregnum(M)+i.$
So we recover the usual result in the standard
graded case.
\end{rmk}

We need two familiar results about $\vregnum(M)$ that follow from results in \cite{MS}.

\begin{lem}\label{L:exact3}
If
\[0 \to M' \to M \to M'' \to 0\] is a short exact sequence of finitely generated
$\ZZ^r$-graded modules, then
\begin{enumerate}
\item $\vregnum(M'') \leq \max\{\vregnum(M), \vregnum(M')-c_{\bf v}\}.$
\item $\vregnum(M) \leq \max\{\vregnum(M'),\vregnum(M'')\}.$
\item $\vregnum(M') \leq \max\{\vregnum(M), \vregnum(M'')+c_{\bf v}\}.$
\end{enumerate}
\end{lem}

\begin{proof}
The result follows from Lemma 7.1 of \cite{MS} and the definition
of $\vregnum(M)$.
\end{proof}

\begin{lem}\label{thm: gens}
Let $M$ be a finitely generated $\ZZ^r$-graded module.
Then the minimal generators of $M^{[{\bf v}]}$ have degrees at most $\vregnum(M)+c_{\bf v}-1.$
\end{lem}

\begin{proof}
The result follows from Theorem 5.4 in \cite{MS} since we have taken
$c_{\bf v} = \lcm(\vdeg (x_i))_{i=1}^n.$
\end{proof}

\begin{proof}[Proof of Theorem \ref{thm: free res}]
We construct a minimal free resolution ${\bf F}_{\bullet}$ satisfying the claim as follows.
By Lemma \ref{thm: gens} we know that the minimal generators of $M$ have
degree at most $\vregnum(M) +c_{\bf v}-1.$  So we have
\[ 0 \to M_1 \to F_0 \to M \to 0\]
and we see that
\begin{eqnarray*}
\vregnum(M_1) &\leq& \max \{ (n-1)c_{\bf v} - \sum \vdeg(x_i) +1  + (\vregnum(M)+c_{\bf v}-1),\\
&&\vregnum(M)+c_{\bf v}\}  \\
&=&  \max \{ nc_{\bf v} - \sum \vdeg(x_i)+ \vregnum(M), \vregnum(M)+c_{\bf v}\}\\
&=& s_{\bf v} + \vregnum(M).
\end{eqnarray*}

Now we can proceed by induction on the projective dimension of $M$. There exists
a short exact sequence
\[ 0 \to M_{i+1} \to F_i \to M_i \to 0\] where
$\vregnum(M_i)\leq is_{\bf v} +\vregnum(M).$ Then $F_i$ is generated by elements of
degree  at most $is_{\bf v} +\vregnum(M)+c_{\bf v}-1.$  Furthermore,
\small
\begin{eqnarray*}
\vregnum(M_{i+1}) &\leq& \max\{(n-1)c_{\bf v} - \sum \vdeg(x_i)+1 +(is_{\bf v} +\vregnum(M)+c_{\bf v}-1),\\
&&is_{\bf v} +\vregnum(M)+c_{\bf v}\}\\
&=& \max\{ nc_{\bf v} - \sum \vdeg(x_i) + is_{\bf v} +\vregnum(M), \\
&&is_{\bf v}+\vregnum(M)+c_{\bf v}\}\\
&=& (i+1)s_{\bf v} +\vregnum(M).
\end{eqnarray*}
\normalsize
We conclude that $F_{i+1}$ is generated by elements of degree at  most \[(i+1)s_{\bf v}
+\vregnum(M)+c_{\bf v}-1\] by Lemma \ref{thm: gens}.
\end{proof}

\begin{cor}\label{cor: minimal_generator}
Let ${\bf F}_{\bullet}$ be a minimal free $\ZZ^r$-graded resolution of $M.$  Assume that the
Hilbert series of $M$ is supported on  $ \bigcup ( {\bf b}_i +Q) $ where ${\bf b}_1,
\ldots, {\bf b}_k\in \ZZ^r.$  Then the minimal
generators of $F_i$ have multidegrees contained in the finite set
\[\D_{i,{\bf v}}(M) :=
\left.\left\{{\bf a} \in \bigcup ({\bf b}_i +Q) ~\right|~ {\bf a}\cdot {\bf v} \leq \vregnum(M)+
is_{\bf v}+c_{\bf v}-1 \right\} .\]
\end{cor}

\subsection{Minimal sets of positive coarsening vectors}

If ${\bf v}$ is a positive coarsening vector,  then
the set $\D_{i,\bf v}(M)$ of Corollary  \ref{cor: minimal_generator}
bounds the multidegrees of the generators of $F_i$, the $i$th syzygy module
of $M$.  In fact, since Corollary \ref{cor: minimal_generator} holds for every positive
coarsening vector ${\bf v},$ we see that we may strengthen the result as follows:
the minimal
generators of $F_i$ have multidegrees contained in
\[\D_i(M) := \bigcap_{{\bf v}} \D_{i,{\bf v}}(M) =\bigcap_{{\bf v}}
\left.\left\{ {\bf a} \in \bigcup ({\bf b}_k +Q) ~\right|~ {\bf a}\cdot {\bf v}
\leq \vregnum(M)+is_{\bf v}+c_{\bf v}-1 \right\} \] where the
intersection is taken over all vectors ${\bf v}$ that are positive
coarsening vectors.

Although there are an infinite number of positive coarsening vectors, the
region $\D_i(M)$ may be computed from only finitely many positive coarsening
vectors ${\bf v}.$

\begin{thm}\label{thm: mspcv}
Let $M$ be a finitely generated $\ZZ^r$-graded $S$-module.
Then there exists a finite set of positive coarsening
vectors $\mathcal{V} = \{{\bf v}_1, \ldots, {\bf v}_n\}$ such that
\[\D_i(M) = \bigcap_{{\bf v}_j \in \mathcal{V}}  \left.\left\{ {\bf a} \in \bigcup ({\bf b}_k +Q)
~\right|~
{\bf a}\cdot {\bf v}_j \leq \vregnum_j(M)+is_{{\bf v}_j}+c_{{\bf v}_j}-1 \right\}.\]
\end{thm}

\begin{proof}
Let ${\bf v}$ be any positive coarsening vector for $M.$  Then $\D_{i, {\bf v}}(M)$
contains only finitely many integral vectors.  We know that $\D_i(M) \subseteq \D_{i, {\bf v}}(M)$
by definition.  Since $\D_{i,{\bf v}}(M)$ is finite, we see that the
set $\D_{i,{\bf v}}(M) - \D_i(M)$ must also be finite.
For each ${\bf p} \in \D_{i,{\bf v}}(M) - \D_i(M)$ there must
be some positive coarsening vector ${\bf v}_{\bf p}$ such that
\[ {\bf p} \cdot {\bf v}_{\bf p} > \vregnum_{\bf p}(M) + is_{{\bf v}_{\bf p}} +c_{{\bf v}_{\bf p}}-1.\]
The set of all ${\bf v}_{\bf p}$ together with ${\bf v}$ is finite and $\D_i(M)$
may be computed from these vectors.
\end{proof}

In light of the previous result, it is natural to introduce the following definition:

\begin{defin} \label{dfn: mspcv}
 Let $M$ be a finitely generated $\ZZ^r$-graded $S$-module.
We shall say a set  $\mathcal{V} = \mathcal{V}(M) = \{{\bf v}_1,\ldots,{\bf v}_n\} \subseteq \ZZ^r$ is
a \emph{minimal set of positive coarsening vectors} for $M$ if
\begin{enumerate}
\item[$(1)$] each ${\bf v}_i$ is a positive
coarsening vector of $S$,
\item[$(2)$]  $\D_i(M) = \bigcap_{{\bf v}_j \in \mathcal{V}} \D_{i,{\bf v}_j}(M)$,
\item[$(3)$] if $\mathcal{V}'$ is a proper subset of $\mathcal{V}$, then
$\D_i(M) \neq \bigcap_{{\bf v} \in \mathcal{V}'} \D_{i, \bf v}(M)$.
\end{enumerate}
\end{defin}

\begin{rmk}As we shall see in \S 4, it appears quite difficult to find a minimal set of positive
coarsening vectors (see Remark \ref{rm: difficult} for more details).  We therefore
propose some natural questions about this minimal set:  is there a
method to identify $\mathcal{V}$ from the $\ZZ^r$-grading?  is the
size of $\mathcal {V}$ an invariant?  is $\mathcal{V}$ unique?
\end{rmk}

\subsection{Further properties of the regularity number}

It is clear from the above discussion that $\vregnum(M)$ gives information
about the $\ZZ$-graded resolution of $M^{[{\bf v}]}$. In the standard graded case, the
process can be reversed, that is, the regularity of $M$ can be read off of the resolution.
As the following
result and example show, we can obtain a lower bound on $\vregnum(M)$ from the
$\ZZ$-graded resolution of $M^{[{\bf v}]}$, but we cannot
determine an upper bound on $\vregnum(M)$ from the invariants of the resolution.

\begin{thm}
Let $M$ be a finitely generated positively $\ZZ^r$-graded $S$-module,
and let ${\bf v}$ be a positive coarsening vector.  If
\begin{equation*}
0 \rightarrow F^{[{\bf v}]}_{\ell} \rightarrow
 F^{[{\bf v}]}_{\ell-1} \rightarrow \cdots
\rightarrow F^{[{\bf v}]}_{1} \rightarrow
 F^{[{\bf v}]}_{0} \rightarrow M^{[{\bf v}]} \rightarrow 0
\end{equation*}
where $F^{[{\bf v}]}_i =  \bigoplus_{j \in \ZZ} S^{[{\bf v}]}(-j)^{\beta_{i,{j}}(M^{[{\bf v}]})}$
is a minimal $\ZZ$-graded resolution of $M^{[{\bf v}]}$, then
\[\max_{i,j}\{j - is_{\bf v} - c_{\bf v} + 1 ~|~ \beta_{i,j}(M^{[{\bf v}]}) \neq 0 \} \leq
\vregnum(M).\]
\end{thm}

\begin{proof}
Suppose $\beta_{i,j}(M^{[{\bf v}]}) \neq 0$.  So, there exists a minimal generator of the
$i$th syzygy module of degree $j$.  By Theorem \ref{thm: free res} we must have
\[j \leq \vregnum(M) + is_{\bf v} + c_{\bf v} -1 \Leftrightarrow
j - is_{\bf v} - c_{\bf v} + 1 \leq \vregnum(M).\]
The above inequality is true for all $i,j$ with $\beta_{i,j}(M^{[{\bf v}]}) \neq 0$, so it
also holds for the maximum such value.
\end{proof}

When ${\bf v} = (1,\ldots,1)$, then $\vregnum(M)$ is simply the Castelnuovo-Mumford regularity
of $M$, and in this case, $\vregnum(M)$ equals
$\max_{i,j}\{j - is_{\bf v} - c_{\bf v} + 1 ~|~ \beta_{i,j}(M^{[{\bf v}]}) \neq 0 \}$.  In fact,
this is often taken as a definition.  In general, $\vregnum(M)$
may be bigger than this lower bound as exhibited in the following example.

\begin{ex}
Let $S = k[x,y]$ with $\deg(x) =(1,0)$ and $\deg(y) = (0,1)$.  Let $M$ be the
free $S$-module $M = S(-1,0) \oplus S(0,-1)$.  Now $M$ has projective dimension equal to
zero, so a resolution of $M$ is simply
\[0 \rightarrow S(-1,0) \oplus S(0,-1) \rightarrow M \rightarrow 0.\]
A positive coarsening vector ${\bf v}$ of $S$ must have the form ${\bf v} = (v_1,v_2) \geq (1,1)$.
Thus, a resolution of $M^{[{\bf v}]}$ as a $\ZZ$-graded module is
\[0 \rightarrow S^{[{\bf v}]}(-v_1) \oplus  S^{[{\bf v}]}(-v_2) \rightarrow M^{[{\bf v}]} \rightarrow 0.\]
Then
\[\max\{j-is_{\bf v} - c_{\bf v} + 1 ~|~ \beta_{i,j}(M^{[{\bf v}]}) \neq 0 \}
= \max\{v_1 - c_{\bf v} + 1, v_2 - c_{\bf v} + 1 \}.\]
On the other hand, by Corollary \ref{C:pdim0} we have
\[\vregnum(M) = \vregnum(S) + \max\{v_1,v_2\} = c_{\bf v} + 1 - v_1 - v_2 + \max\{v_1,v_2\}.\]
Suppose that ${\bf v} = (v_1,v_2)$ has been picked so that
gcd$(v_1,v_2) = 1$ and $v_1>v_2 > 1$ (such pairs exist, e.g.,
${\bf v} = (5,3)$).
Then $c_{\bf v} = v_1v_2$ and
\begin{eqnarray*}
\max\{j-is_{\bf v} - c_{\bf v} + 1 ~|~ \beta_{i,j}(M^{[{\bf v}]}) \neq 0 \} & = &
v_1 - v_1v_2 + 1 = v_1(1-v_2)+1\\
&<& v_1v_2 - v_2 + 1 = \vregnum(M).
\end{eqnarray*}
\end{ex}
Note that we may also define the familiar $a$-invariants of $M^{[{\bf v}]}$ in this context as well:

\begin{defin}
Let $M$ be a finitely generated positively $\ZZ^r$-graded $S$-module and let ${\bf v}$ be a positive
coarsening vector.  For each $i \geq 0$, let
\[
 a^i(M^{[{\bf v}]}) := \max \{p \mid H^i_{\mathfrak{m}}(M^{[{\bf v}]})_p \neq 0 \}
\]
where $a^i(M^{[{\bf v}]}) = -\infty$ if $H^i_{\mathfrak{m}}(M^{[{\bf v}]}) = 0$.
\end{defin}
\begin{rmk}
Benson \cite{B} defines a notion of regularity for a module $M$ over
a graded commutative
ring $H$ with nonstandard $\ZZ$-grading to be $\max \{a^i(M) +i \}.$  However,
to link a notion of regularity defined in terms of local cohomology with
degrees of generators as in Lemma \ref{thm: gens}, the standard proof
proceeds by induction on the dimension of $M$ via a ring element $x$ that is
a nonzero divisor (or almost a nonzerodivisor) on $M$  (see \cite{B},
\cite{BS}, and \cite{MS}.)  Since $c_{\bf v}$ is the least
common multiple of the degrees of the variables $x_i,$ we know we can find
an element that is almost a nonzerodivisor in degree $c_{\bf v}$ from Proposition 3.1 of \cite{MS}.  This is why we require vanishings in
local cohomology shifted by $c_{\bf v}.$  (See also Theorem 5.4 of \cite{MS}.)
Benson links his notion of regularity to free resolutions by passing to
a Noether normalization $R$ of $H$ generated by elements which form a
 filter-regular (almost regular) sequence on $M.$  Under suitable
hypotheses there are nice links between Benson's regularity and free
resolutions of $M$ over $R.$  However, since we are interested in resolutions of $M$ over $S,$ we have chosen to
follow along the lines of \cite{MS}.
\end{rmk}

The quantity $\vregnum(M)$ interacts with the $a$-invariants as follows:

\begin{thm}\label{lem:a-inv}
Keeping notation and hypotheses as above, \[\vregnum(M) = \max_i \{a^i(M^{[{\bf v}]}) -c_{\bf v}(1-i)+1 \}.\]
\end{thm}

\begin{proof}
First, we show that $a^i(M^{[{\bf v}]}) -c_{\bf v}(1-i) +1\leq \vregnum(M)$ for all $i.$
If $r = \vregnum(M),$ then $H^i_{\mathfrak{m}}(M)_p = 0$
for all $p \in r + j + c_{\bf v}(1-i) + \NN c_{\bf v},$ for all $j \in \NN$  and in
particular, for all $p \geq r+c_{\bf v}(1-i).$  Therefore,
\[a^i(M^{[{\bf v}]}) \leq r+c_{\bf v}(1-i)-1\] or $a^i(M^{[{\bf v}]}) -c_{\bf v}(1-i) +1\leq \vregnum(M)$
for all $i.$

Now we show that $\vregnum(M) \leq \max_i \{a^i(M^{[{\bf v}]}) -c_{\bf v}(1-i) +1 \}.$
Suppose that  $p \geq \max_i \{a^i( M^{[{\bf v}]})-c_{\bf v}(1-i) +1\}.$ It suffices to show that
$p \in \vreg(M),$ which holds if and only if
$H^i_{\mathfrak{m}}(M)_q = 0$ for all $i\geq 0$ and for all
\[q \in p + c_{\bf v}(1-i) + \NN c_{\bf v}.\] Since $q > a^i(M^{[{\bf v}]})$ we have
  $\vregnum(M) \leq \max_i \{a^i(M^{[{\bf v}]}) -c_{\bf v}(1-i)+1 \}.$
\end{proof}


\section{Positive coarsening vectors and their scalar multiples}

In this section we show that if ${\bf v}$ is a positive
coarsening vector, then
no new information on the degrees of the generators is obtained by
using the positive coarsening vector $d{\bf v}$ with $d \in \NN_{>0}$.  As a consequence,
if ${\bf v}  = (v_1,\ldots,v_r)$ and $\ell = \operatorname{gcd}(v_1,\ldots,v_r)$, then ${\bf v}$ can be replaced
with ${\bf v}' = (v_1/\ell, \ldots,v_r/\ell)$.

If $M$ is a finitely generated $\ZZ^r$-graded $S$-module,
and if ${\bf v}$ is a positive coarsening vector, then
\[
M^{[{\bf v}]} := \bigoplus_{m \in \ZZ} \left(\bigoplus_{{\bf a}\cdot{\bf v} =m} M_{\bf a} \right)
~~\mbox{and}~~
M^{[d{\bf v}]} := \bigoplus_{m \in \ZZ} \left(\bigoplus_{{\bf a}\cdot d{\bf v} =m} M_{\bf a} \right)
\]
As a $\ZZ$-graded module, the degree $j$ piece of $M^{[d{\bf v}]}$ is given by
\[
(M^{[d{\bf v}]})_j
= \left\{
\begin{array}{ll}
0 & \mbox{if $j \not\equiv 0 \pmod{d}$} \\
M^{[{\bf v}]}_{\frac{j}{d}} & \mbox{if $j \equiv 0 \pmod{d}$} \\
\end{array}
\right..
\]

\begin{lem} \label{lem: scalar vanishing}  Let $M$ be a finitely generated $\ZZ^r$-graded $S$-module.
For all $j \in \ZZ$
\[H^i_{\mathfrak{m}}(M^{[d{\bf v}]})_j =
\left\{
\begin{array}{ll}
0 & \mbox{if $j \not\equiv 0 \pmod{d}$} \\
H^i_{\mathfrak{m}}(M^{[{\bf v}]})_{\frac{j}{d}} & \mbox{if $j \equiv 0 \pmod{d}$} \\
\end{array}
\right..
\]
\end{lem}

\begin{proof}
The isomorphism  $\phi: M^{[{\bf v}]} \rightarrow M^{[d{\bf v}]}$ where $\phi(f) = f$ satisfies
$\phi((M^{[{\bf v}]})_j) = (M^{[d{\bf v}]})_{dj}$ and gives rise to the
corresponding isomorphism of local cohomology modules.
\end{proof}

\begin{thm}
 Let $M$ be a finitely generated $\ZZ^r$-graded $S$-module.
If $p \in \reg_{\bf v}(M^{[{\bf v}]})$, then
\[dp - \ell \in \reg_{d{\bf v}}(M^{[d{\bf v}]}) ~~\mbox{for $\ell = 0,\ldots,d-1$}.\]
\end{thm}

\begin{proof}
We will consider the case $\ell =0$ and $\ell \neq 0$ separately.
Suppose $\ell = 0$.  We want to show that $dp \in \reg_{d{\bf v}}(M^{[d{\bf v}]})$.
Fix an $i \geq 0$ and let
\[ q \in dp + \NN dc_{\bf v}[1-i].\]
So, $q = dp + dc_{\bf v}(1-i+j)$ for some $j \in \NN$.  Since $q \equiv 0 \pmod{d}$ for
all $i,j$, we have by Lemma \ref{lem: scalar vanishing}
\[H^i_{\mathfrak{m}}(M^{[d{\bf v}]})_q = H^i_{\mathfrak{m}}(M^{[{\bf v}]})_{\frac{q}{d}}
= H^i_{\mathfrak{m}}(M^{[{\bf v}]})_{p + c_{\bf v}(1-i+j)}.\]
But $p \in \reg_{\bf v}(M^{[\bf v]})$ so $H^i_{\mathfrak{m}}(M^{[{\bf v}]})_{p + c_{\bf v}(1-i+j)} = 0$.
Thus $dp \in \reg_{d{\bf v}}(M^{[d{\bf v}]})$.

Suppose now that  $\ell \neq 0$ and fix an $i \geq 0$.  Let
\[ q \in (dp - \ell) + \NN dc_{\bf v}[1-i].\]
But then $q = dp - \ell + dc_{\bf v}(1-i+j)$ for some $j \in \NN$, and hence $q \not\equiv 0 \pmod{d}$.
By Lemma \ref{lem: scalar vanishing} we get
$H^i_{\mathfrak{m}}(M^{[d{\bf v}]})_q = 0,$ which implies that
$dp - \ell \in \reg_{d{\bf v}}(M^{[d{\bf v}]})$.
\end{proof}

\begin{cor}\label{C:dv}
 Let $M$ be a finitely generated $\ZZ^r$-graded $S$-module.  Then
\[\regnum_{d{\bf v}}(M^{[d{\bf v}]}) = d\vregnum(M^{[{\bf v}]}) - d + 1.\]
\end{cor}

\begin{proof}
Let $r = \vregnum(M^{[{\bf v}]})$, and thus $r' \in \reg_{\bf v}(M^{[{\bf v}]})$ for
all $r' \geq r$.  The previous theorem implies that for all $s \geq dr - d + 1$,
$s \in \reg_{d\bf v}(M^{[{d\bf v}]})$.  So
\[dr - d + 1 \geq \regnum_{d{\bf v}}(M^{[d{\bf v}]}).\]
On the other hand, we know that $r - 1 \not\in  \reg_{\bf v}(M^{[{\bf v}]})$.  If
$d(r-1) \in \reg_{d{\bf v}}(M^{[d{\bf v}]})$ this would imply that
\[H^i_{\mathfrak{m}}(M^{[d{\bf v}]})_q = 0 ~~\forall i\geq 0 ~~\forall q \in d(r-1) + \NN
dc_{\bf v}[1-i].\]
But by Lemma \ref{lem: scalar vanishing}, this would imply
\[H^i_{\mathfrak{m}}(M^{[{\bf v}]})_{q'} = 0 ~~\forall i\geq 0 ~~ \forall q' \in (r-1) + \NN c[1-i]\]
contradicting the fact that $r -1 \not\in \reg_{\bf v}(M^{[{\bf v}]})$.
Combining the inequality
\[\regnum_{d{\bf v}}(M^{[d{\bf v}]}) > d(r-1) = dr -d \]
with our previous inequality gives the
desired conclusion.
\end{proof}

We can now show that if we have
 a bound on the multidegrees of the generators using ${\bf v}$,
we will not obtain new bounds on the multidegrees if we use a scalar multiple of ${\bf v}$.

\begin{thm}\label{T:dv}
Let ${\bf v}$ be a positive coarsening vector and let $d \in \NN_{> 0}$. Then
\[\mathcal{D}_{i, d{\bf v}}(M) = \mathcal{D}_{i,{\bf v}}(M).\]
\end{thm}
\begin{proof}
A simple calculation will show that $c_{d \bf v} = dc_{\bf v}$ and $s_{d{\bf v}} = d s_{\bf v}$.
Suppose that ${\bf a} \in \mathcal{D}_{i,d {\bf v}}(M)$.  So
\begin{eqnarray*}
{\bf a}\cdot d{\bf v} &\leq& \regnum_{d{\bf v}}(M^{[d{\bf v}]}) + is_{d{\bf v}} + c_{d\bf v} -1 \\
&=&  d\regnum_{\bf v}(M^{[\bf v]}) - d + 1 + ids_{\bf v} + dc_{\bf v} -  1.
\end{eqnarray*}
Dividing both sides of the inequality by $d$ gives
\[{\bf a}\cdot {\bf v} \leq  \regnum_{\bf v}(M^{[\bf v]}) + is_{\bf v} + c_{\bf v} -  1,\]
that is, ${\bf a} \in \mathcal{D}_{i,{\bf v}}(M)$.  By reversing the argument, we
also can show that $\mathcal{D}_{i,{\bf v}}(M) \subseteq \mathcal{D}_{i,d{\bf v}}(M)$.\end{proof}

\begin{rmk}\label{rm: difficult}
The result above tells us that a minimal set of positive coarsening vectors will not contain parallel vectors.  However, to proceed further, we need a mechanism to compare the sets $\D_{i,{\bf v}}(M)$ and $\D_{i,{\bf w}}(M)$
for two arbitrary positive coarsening vectors ${\bf v}$ and ${\bf w}$.
 The primary obstacle in
describing a comparison is relating the vanishing of $H_{\mathfrak{m}}^i(M^{[{\bf v}]})$
with that of $H_{\mathfrak{m}}^i(M^{[{\bf w}]})$, as
we did in Lemma \ref{lem: scalar vanishing} for ${\bf v}$ and $d{\bf v}$.
In particular, we want a relationship like Corollary \ref{C:dv}
for $\vregnum(M)$ and $\operatorname{reg-num}_{{\bf w}}(M)$.
\end{rmk}


\section{Connections with other notions of regularity}\label{sec: B-reg}

In this section we relate $\vregnum(M)$ to two recent notions of multigraded regularity.
In particular, we show how to use $\vregnum(M)$ to determine bounds on the
multigraded regularity of $M$ as defined in \cite{MS}.  As well, we relate
$\vregnum(M)$ to the resolution regularity vector of $M$ as defined by \cite{ACD,SV}.

\subsection{$B$-regularity}\label{section: B-reg}
We assume that our multigraded ring $S$ is the homogeneous coordinate ring
 of a smooth toric variety $X$ of dimension $d = n-r$ whose Chow group
$A_{d-1}(X) \cong \ZZ^r$ for
some $r.$  The ring $S$ comes equipped with a square-free
``irrelevant'' monomial ideal $B.$  The group in which
the degrees of $S$ lie is the Chow group $A_{d-1}(X).$

Alternatively, we may work in the more general setting described in \S 2 of
\cite{MS} by imposing the additional condition that $S$ is positively multigraded
by $\ZZ^r$ and that $B$ is chosen in a way that is compatible with the grading.

Assume that an isomorphism $A_{d-1}(X) \cong \ZZ^r$ and a set
$\mathcal{C} =
\{{\bf c}_1, \ldots, {\bf c}_{\ell}\} \subset \ZZ^r$ have been chosen.
We let $\NN\mathcal{C}$ denote the subsemigroup of $\ZZ^r$ generated by $\mathcal{C}$.

\begin{defin}[Definition 4.1 in \cite{MS}]
Let $M$ be a finitely generated
$\ZZ^r$-graded $S$-module.  If ${\bf m} \in \ZZ^r$, then ${\bf m}$ is in the \emph{$B$-regularity} of $M$ from
level $i$ with respect to $\mathcal{C}$ if
\[H^j_B(M)_{\bf p} = 0 \ ~\mbox{for all}~ j \ge i ~\mbox{and } {\bf p}\in {\bf m}
+\NN \mathcal{C}[1-j]\] where
\[\NN \mathcal{C}[j] = \bigcup_{{\bf w} \in \NN^{\ell}, \sum w_i = |j| }
\sign(j) (w_1{\bf c}_1+\cdots+w_{\ell}{\bf c}_{\ell}) +\NN \mathcal{C}. \]
The set of all such ${\bf m}$ is denoted $\reg_B^i(M)$ and $\reg_B(M) = \reg_B^0(M).$
\end{defin}

The $B$-regularity of the polynomial ring $S$
can always be computed topologically via Proposition 3.2 of \cite{MS}.
Once the sets $\reg_B^i(S)$ are known, we can compute
a lower bound on $\reg_B(M)$ for any finitely generated $\ZZ^r$-graded
$S$-module from the coarsely graded minimal free resolution of
$M$, and hence, from $\vregnum(M).$
\begin{thm}\label{thm: res bound}
Let $M$ be a finitely generated $\ZZ^r$-graded $S$-module.  If ${\bf v}$ is a
positive coarsening vector, and if
\[{\bf F}_{\bullet}: \cdots \to F_i  \to \cdots \to F_0 \to M^{[{\bf v}]} \to 0 \]
is a minimal free $\ZZ$-graded resolution of $M^{[{\bf v}]}$ with $F_i = \bigoplus_j
S^{[{\bf v}]}(-d_{i,j}),$
then
\[ D_0 \cap
\left( \bigcup_{1\leq j \leq \ell} -{\bf c}_j + 
\bigcup_{\phi:[d+1] \to [\ell]} \left(
\bigcap_{1 \leq i\leq d+1} (- {\bf c}_{\phi(2)} - \cdots - {\bf c}_{\phi(i)} + D_i) \right)
\right) \subseteq \reg_B(M)\]
where \[D_i = \bigcap_{{\bf a}\cdot {\bf v} = d_{i,j} \ \mathrm{for} \ \mathrm{some} \ j}
\bigl({\bf a} + \reg_B^i(S)\bigr).\]
\end{thm}

\begin{proof}
 We know that $M$ has a minimal free $\ZZ^r$-graded resolution
\[0 \to \cdots \to G_i \to \cdots \to G_0 \to M \to 0\]
where $G_i = \bigoplus_k S(-{\bf a}_{i,k}).$  Coarsening this
resolution with the vector ${\bf v}$ gives a minimal free
$\ZZ$-graded resolution of $M^{[{\bf v}]}$ which must be
isomorphic to ${\bf F}_{\bullet}$. One can easily see that Lemma
A.1 and Theorem A.2 in \cite{SV} (variants of Lemma 7.1 and
Theorem 7.2 in \cite{MS}) hold in this slightly more general
situation, so we know that
{\footnotesize
\[ D_0 \cap \left(
\bigcup_{1 \leq j \leq \ell} -{\bf c}_j +   \bigcup_{\phi:[d+1] \to [\ell]}
\left( \bigcap_{1 \leq i\leq d+1} (- {\bf c}_{\phi(2)} - \cdots - {\bf c}_{\phi(i)} + \reg_B^i(G_i))\right)\right) \subseteq \reg_B(M).\]}
So it is enough to show that $D_i \subseteq \reg_B^i(G_i).$  But
\[ \reg_B^i(G_i) = \bigcap_k \bigl( {\bf a}_{i,k} + \reg_B^i(S) \bigr),\]
which contains $D_i$
since for each $k$ such that ${\bf a}_{i,k} \neq 0,$ there exists a $d_{i,j}$
in the $\ZZ$-graded resolution of $M^{[{\bf v}]}$ with
${\bf a}_{i,k} \cdot {\bf v} = d_{i,j}.$
\end{proof}

\begin{cor}
Under the same hypotheses as Theorem \ref{thm: res bound},
{\footnotesize
\[ \bigcup_{\phi:[d+1] \to [\ell]} D_0 \cap
\left( \bigcup_{1 \leq j \leq \ell} -{\bf c}_j +\bigcup_{\phi:[d+1] \to [\ell]}
\left( \bigcap_{1 \leq i\leq d+1} (- {\bf c}_{\phi(2)} - \cdots - {\bf c}_{\phi(i)} + D_i) \right)\right) \subseteq \reg_B(M)\]}
where \[D_i = \bigcap_{{\bf a}\cdot {\bf v} \leq \vregnum(M)+is_{\bf v} +{\bf c}_{\bf v}-1}
\bigl({\bf a} + \reg_B^i(S)\bigr).\]
\end{cor}
\begin{proof}
Keeping the setup of the proof of Theorem \ref{thm: res bound},
\[\{{\bf a} \mid {\bf a}\cdot {\bf v} \leq \vregnum(M)+is_{\bf v} +c_{\bf v}-1\} \supseteq\{{\bf a} \mid {\bf a} \cdot {\bf v} = d_{i,j} \ \mathrm{for} \ \mathrm{some} \ j\}. \]
\end{proof}

\subsection{Resolution regularity vector}
We shall now assume that
$S$ is
the standard multigraded homogeneous coordinate ring of $\mathbb{P}^{n_1} \times \cdots \times \mathbb{P}^{n_r}$.
Generalizing the bigraded case found in \cite{ACD}, the first two authors introduced the
following definition in \cite{SV}:

\begin{defin}\label{def: res-reg}
Let $M$ be a finitely generated $\ZZ^r$-graded $S$-module.
For $\ell = 1,\ldots,r$, let
\[ d_{\ell} := \max \{a_{\ell}~|~\tor_i^S(M, k)_{(a_1, \ldots, a_{\ell}+i, \ldots, a_r)} \neq 0  \}\]
for some $i$ and for some $a_1, \ldots, a_{\ell-1}, a_{\ell+1}, \ldots, a_r.$  We
call $\underline{r}(M) := (d_1, \ldots, d_r)$ the \emph{resolution regularity vector} of $M.$
\end{defin}

It follows that if $\underline{r}(M) = (d_1,\ldots,d_r)$ and if $F$ if is a minimal
generator of the $i$th syzygy module of $M$, then the $\ell$th coordinate of $\deg F$ is
bounded above $d_{\ell} + i$.

\begin{thm}\label{thm: connections}
Let $M$ be a finitely generated $\ZZ^r$-graded $S$-module whose Hilbert series is
supported on $\NN^r,$ and let
\[ {\bf F}_{\bullet}: ~\cdots \to F_i \to \cdots \to F_0 \to M \to 0\]
be a minimal free $\ZZ^r$-graded resolution of $M.$
Let $d = \vregnum(M)$ when ${\bf v} = (1,\ldots,1)$.
\begin{enumerate}
\item The free module $F_i$ is generated by elements whose degrees are contained in the set of all
${\bf a} \in \NN^r$ satisfying $a_1+\cdots + a_r \leq d+i.$
\item  Set $m = \min\{\pdim(M), n_1+\cdots + n_r+1\}.$  If $m >0,$ then
\[\underline{r}(M) \leq (d,\ldots, d) ~\mbox{and}~
(d+m,\ldots, d+m)+\NN^r[-(m-1)] \subseteq \reg_B(M).\]
\end{enumerate}
\end{thm}

\begin{proof}
Statement $(1)$ follows from Corollary \ref{cor: minimal_generator}.  For
$(2)$, the first part is clear, and the second follows from Corollary 2.3 of \cite{SV}.
\end{proof}

\section{Illustrative Examples}

We end this paper with some examples that illustrate the strengths
and weaknesses of using coarsening vectors to find bounds on the degrees of
the multigraded generators of the syzygies in the $\ZZ^r$-graded free
resolutions of finitely generated $\ZZ^r$-graded $S$-modules.

\subsection{Points in multi-projective spaces}
For our first set of examples, we consider the ideals of sets of points in
$\PP^{n_1} \times \cdots \times \PP^{n_r}$.
The benefit of studying this class of ideals is that we already have some information on the regularity of these
ideals (see \cite{HV,MS,SV} for more on this).

Let $X$ be a finite set of reduced points in $\PP^{n_1} \times \cdots \times \PP^{n_r}$, that
is, the defining ideal of $X$ equals its radical.
The defining ideal of $X$, denoted $I_X$, is an $\NN^r$-graded homogeneous
ideal of the $\NN^r$-graded polynomial ring
$S = k[x_{1,0},\ldots,x_{1,n_1},\ldots,x_{r,0},\ldots,x_{r,n_r}]$
where $\deg (x_{i,j}) = e_i$, the $i$th standard basis vector
of $\NN^r$.  The {\it Hilbert function} of $X$ is the numerical function
$H_X:\NN^r \rightarrow \NN$ defined by ${\bf i} \mapsto \dim_k (S/I_X)_{\bf i}$.

As  shown in \cite{MS}, the $B$-regularity of $S/I_X$ is linked
to $H_X$.  Specifically,

\begin{thm}\label{MS-points}
Let $X$ be a set of reduced points in $\PP^{n_1} \times \cdots \times \PP^{n_r}$.  Then
\[\reg_B(S/I_X) = \{{\bf i} \in \NN^r ~|~ H_X({\bf i}) = \deg X = |X|\}.\]
\end{thm}

The resolution regularity vector of $S/I_X$ can also be read off of the Hilbert function
of $S/I_X$.  The following theorem is a special case of Theorem 4.2 of \cite{SV}.

\begin{thm}\label{SV-points}
Let $X$ be a set of reduced points in $\PP^{n_1} \times \cdots \times \PP^{n_r}$.  Then
\[\underline{r}(S/I_X) = (t_1,\ldots,t_r)\]
where $t_i = \min\{t \in \NN ~|~ H_X(te_i) = |\pi_i(X)|\}$ and $\pi_i:\PP^{n_1} \times
\cdots \times \PP^{n_r} \rightarrow \PP^{n_i}$ is the $i$th projection morphism.
\end{thm}

The following two examples compare these notions of regularity with the notion
developed in this paper.

\begin{ex}
In $\PP^1 \times \PP^1$, let $P_{ij} = [1:i] \times [1:j]$.  Consider the
set of points $X = \{P_{11},P_{12},P_{13},P_{14},P_{21},P_{31},P_{41},P_{51}\}$.  The
Hilbert function of $X$ is
\[
\left[
\begin{matrix}
\vdots & \vdots & \vdots & \vdots & \vdots & \vdots& \\
4 & 5 & 6 & 7 & 8 &8&\cdots \\
4 & 5 & 6 & 7 & 8 &8&\cdots \\
3 & 4 & 5 & 6 & 7 & 7&\cdots \\
2 & 3 & 4 & 5 & 6 & 6&\cdots\\
1 & 2 & 3 & 4 & 5 & 5& \cdots
\end{matrix}
\right],
\]
where the number in the location $(i,j)$ corresponds to the value of $H_{X}(i,j)$.  By Theorem
\ref{MS-points}, $\reg_B(S/I_X) = (4,3) + \NN^2.$
The degrees of the minimal $i$th syzygies then belong to the unbounded set
\[\mathcal{A}_i  =\NN^2 \backslash \bigcup_{m,n \geq 0, m+n=i} ((4+m,3+n)  + \NN^2).\]

Theorem ~\ref{SV-points} tells us that
 $\underline{r}(S/I_X) = (4,3)$ since $4 = \min\{t ~|~ H_{X}(t,0) = |\pi_1(X)| = 5\}
$ and $3=\min\{t ~|~ H_{X}(0,t) = |\pi_2(X)| = 4\}$.  The resolution regularity
vector tells us that the degree of a minimal $i$th syzygy belongs to the finite set
\[\mathcal{B}_i = \{(m,n) \in \NN^2 ~|~ (0,0) \leq (m,n) \leq (4+i,3+i)\}. \]
Observe that $\mathcal{A}_i \not\subset\mathcal{B}_i$ and $\mathcal{B}_i \not\subset \mathcal{A}_i$
since if $i=2$, $(6,5) \in \mathcal{B}_2$ but is not in $\mathcal{A}_2$, while $(7,0) \in \mathcal{A}_2$
but not in $\mathcal{B}_2$.

Let us now consider the coarsening vector ${\bf v} = (1,1)$. Under this
coarsening vector  $\vregnum(S/I_X)$ agrees with the standard notion of regularity of
a $\ZZ$-graded module.  Thus $\vregnum(S/I_X) = 4$ since the
$\ZZ$-graded resolution of $(S/I_X)^{[{\bf v}]}$ has the form
\[0 \rightarrow S^{[{\bf v}]}(-5)\oplus S^{[{\bf v}]}(-6) \rightarrow S^{[{\bf v}]}(-5)
\oplus S^{[{\bf v}]}(-4) \oplus S^{[{\bf v}]}(-2)
\rightarrow S^{[{\bf v}]} \rightarrow
(S/I_X)^{[{\bf v}]} \rightarrow 0\]
The generators of the $i$th syzygy module of $S/I_X$ must therefore have degrees belonging to the set
\[\mathcal{C}_i = \D_{i,{\bf v}}(S/I_X) = \{(m,n) \in \NN^2 ~|~ m+n \leq 4 + i \}.\]
Note that $\mathcal{B}_i \not\subset \mathcal{C}_i$ and $\mathcal{C}_i \not\subset \mathcal{B}_i$
since $(0,5) \in \mathcal{C}_1$ but $(0,5) \not\in \mathcal{B}_1$, while $(4,5) \in \mathcal{B}_1$
but not in $\mathcal{C}_1$.

On the other hand, $\mathcal{C}_i \subseteq \mathcal{A}_i$ for all $i$.  So, for this example, the
coarse grading yields new information on degrees of the minimal syzygies when compared to
the $B$-regularity of \cite{MS}.
\end{ex}

The previous example showed that a coarse grading can reveal new information about
the degrees of the minimal syzygies.  However, this may not always be the case,
as Example \ref{ex: generic points} will show.  We first recall some relevant
results about points in generic position.

A set of points $X \subseteq \PP^{n_1} \cdots \times \PP^{n_r}$ is
said to be in \emph{generic position} if $H_X({\bf i}) = \min\{\dim_{k} S_{\bf i},|X|\}$
for all ${\bf i} \in \NN^r.$  The regularity of $S/I_X$ as a $\ZZ$-graded ring
was computed in \cite{HV} when $X$ is in generic position:

\begin{thm}\label{HV-points}
Let $X$ be a set of reduced points in $\PP^{n_1} \times \cdots \times \PP^{n_r}$ in
generic position.  Then
\[\reg(S/I_X) = \max\{d_1,\ldots,d_r\}\]
where $d_i = \min\left\{d ~\left|~ \binom{d+n_i}{d} \geq |X| \right\}\right.$
for $i =1,\ldots,r$.
\end{thm}

\begin{ex}\label{ex: generic points}
Let  $X \subseteq \PP^{n_1} \cdots \times \PP^{n_r}$ be
a set of points in generic position.  By Theorem \ref{MS-points},
\[
\reg_B(S/I_X) = \left\{{\bf i}\in \NN^r ~\left|~ \binom{n_1+i_1}{n_1}\cdots\binom{n_r+i_r}{n_r} \geq |X| \right\}
\right..\]
The degrees of the minimal $i$th syzygies must therefore be elements of
\[\mathcal{A}_i = \NN^r \backslash \bigcup_{
\footnotesize
i_1,\ldots,i_r \geq 0, ~~
i_1+\cdots+i_r = i
\normalsize}
(\reg_B(S/I_X) + (i_1,\ldots,i_r))+ \NN^r.\]
For points in generic position, the sets $\mathcal{A}_i$ are finite.

We compare this to the information we obtained from the coarsely graded
ring.  We use ${\bf v} = (1,\ldots,1)$ as our positive coarsening vector.
Then $\vregnum(S/I_X) = \max\{d_1,\ldots,d_r\}$ where $d_i = \min\{d ~|~ \binom{d+n_i}{d} \geq |X|\}$
by Theorem \ref{HV-points}.  Thus
\[\mathcal{C}_i = \D_{i,{\bf v}}(S/I_X) = \{\ba \in \NN^r  ~|~ \ba\cdot {\bf v} \leq \vregnum(R/I_X) + i\}. \]
But for all $i$, $\mathcal{A}_i \subseteq \mathcal{C}_i$,
so no new information about the multidegrees of $i$th syzygies is
obtained by using the positive coarsening vector ${\bf v}
=(1,\ldots,1)$.
\end{ex}

\begin{ex}\label{ex3}
We find a minimal set of positive coarsening vectors for $X =\{[1:0] \times [1:0],[1:0] \times [0:1],[0:1] \times [1:0], [0:1] \times [0:1]\},$ the set of points of Example 1.1.  The defining ideal $I_X = \langle x_0x_1,y_0y_1 \rangle$ is the Stanley-Reisner ideal of the simplicial complex $\Delta$ equal to the 4-cycle with edges $\{x_0, y_0 \},\{x_0, y_1 \},\{x_1, y_0 \},\{x_1, y_1 \}$ (where we label the vertices with the variables from the ring).

Because $I_X$ is a complete intersection, $S/I_X$ is a
Cohen-Macaulay ring of dimension 2.  Thus,
$H^i_{\mathfrak{m}}(S/I_X) =0$ for all $i \neq 2$.  When $i = 2$,
we can compute the fine $\NN^4$-graded Hilbert series of
$H^2_{\mathfrak{m}}(S/I_X)$ using Hochster's theorem (cf. Theorem
13.13 \cite{Mi-S}). In particular, the Hilbert series
$H(H^2_{\mathfrak{m}}(S/I_X);x_0,x_1,y_0,y_1)$ is

\begin{eqnarray*}
& =& 1+\frac{x_0^{-1}}{1-x_0^{-1}}+\frac{x_1^{-1}}{1-x_1^{-1}}+\frac{y_0^{-1}}{1-y_0^{-1}}+\frac{y_1^{-1}}{1-y_1^{-1}}
+\frac{x_0^{-1}}{1-x_0^{-1}}\cdot\frac{y_0^{-1}}{1-y_0^{-1}} \\
&&+\frac{x_0^{-1}}{1-x_0^{-1}}\cdot\frac{y_1^{-1}}{1-y_1^{-1}}
+\frac{x_1^{-1}}{1-x_1^{-1}}\cdot\frac{y_0^{-1}}{1-y_0^{-1}}
+\frac{x_1^{-1}}{1-x_1^{-1}}\cdot\frac{y_1^{-1}}{1-y_1^{-1}}\\
&=& 1+\left(\frac{1}{x_0} + \frac{1}{x_0^2} + \cdots\right) + \left(\frac{1}{x_1} + \frac{1}{x_1^2} + \cdots\right)+
\left(\frac{1}{y_0} + \frac{1}{y_0^2} + \cdots\right)\\
&&+\left(\frac{1}{y_1} + \frac{1}{y_1^2} + \cdots\right)+ \left(\frac{1}{x_0y_0} + \frac{1}{x_0^2y_0} + \cdots\right)
+ \left(\frac{1}{x_0y_1} + \frac{1}{x_0^2y_1} + \cdots\right) \\
&&+ \left(\frac{1}{x_1y_0} + \frac{1}{x_1y_0^2} +
\cdots\right)+\left(\frac{1}{x_1y_1} + \frac{1}{x_1y_1^2} +
\cdots\right).
\end{eqnarray*}
\normalsize

Let  ${\bf v} = (a,b) \in \NN^2$ be a positive coarsening vector for $S/I_X.$
This implies that $a, b \geq 1.$  We can also assume that $\gcd(a,b) = 1$ (from \S 4) and that $a \leq b.$
Since $\deg(x_i) = a$ and $\deg(y_i) = b,$ the $\ZZ$-graded Hilbert series of $H_{\mathfrak{m}}^2((S/I_X)^{[{\bf v}]})$ is
\begin{eqnarray*}
H(H^2_{\mathfrak{m}}((S/I_X)^{[{\bf v}]});t) &=& 1+2(t^{-a}+t^{-2a} + \cdots) + 2(t^{-b}+t^{-2b}+\cdots)\\
&& + 4(t^{-a-b}+t^{-2a-b}+\cdots).
\end{eqnarray*}
Therefore, we see that $a^2((S/I_X)^{[{\bf v}]}) = 0$
and $a^i((S/I_X)^{[{\bf v}]}) = -\infty$ for all $i \neq 2.$ By Theorem \ref{lem:a-inv}
\[\vregnum(S/I_X) = 0-ab(1-2)+1= ab+1.\]
since $c_{\bf v} = ab$.

We can now show that for any positive coarsening vector ${\bf v} = (a,b)$ and for each $i \geq 0$,
\[\D_{i,{\bf 1}}(S/I_X) \subseteq \D_{i,{\bf v}}(S/I_X)\]
where {\bf 1} $=(1,1)$, and consequently, $\mathcal{V} = \{{\bf 1} \}$ is a minimal set of positive
coarsening vectors for $S/I_X$.  Without loss of generality, we can assume gcd$(a,b) =1$ and
$a \leq b$.
If $(x,y) \in \D_{i,{\bf 1}}(S/I_X) \subset \NN^2,$
then $x + y \leq \regnum_{\bf 1}(S/I_X) + is_{\bf 1} + c_{\bf 1} - 1$.  But since
$s_{\bf 1} = c_{\bf 1} = 1$ and $\regnum_{\bf 1}(S/I_X) = 2$ (because $I_X$ is a complete intersection),
we deduce that $x + y \leq 2 +i.$
Thus
\begin{eqnarray*}
(x,y)\cdot(a,b) = ax+ by &\leq & bx + by \leq 2b + bi\\
& \leq & 2ab + abi = (ab+1) + abi + ab -1
\end{eqnarray*}
Because $c_{\bf v} = ab$, $s_{\bf v} \geq c_{\bf v}$, and since $\vregnum(S/I_X) = ab+1$, we have
\[(x,y)\cdot(a,b) \leq \vregnum(S/I_X) +is_{\bf v} + c_{\bf v}  -1,\]
that is, $(x,y) \in \D_{i,{\bf v}}(S/I_X)$, as desired.
\end{ex}

\begin{rmk}
The previous examples suggest that if one is interested in computing bounds on the degrees of
the minimal syzygies,  one might wish to use the notions of $B$-regularity,
the resolution regularity vector, and the regularity number of a
positive coarsening vector together, as one approach is not always better than another.
\end{rmk}

\subsection{Hirzebruch Surface}
Consider the case where $S=k[x_1, x_2, x_3,x_4]$ is the coordinate
ring of the Hirzebruch surface $\FF_s,$ $s \geq 1.$ The ring $S$ is
multigraded by $A=\ZZ^2$ with
$\deg(x_1)= (1,0), \deg(x_2)= (-s,1), \deg(x_3)= (1,0),$ and $\deg(x_4)= (0,1).$

The solid dots in Figure (i) below represent the degrees of elements of $S.$
The positive coarsening vectors of $S$ are the interior points (represented
as solid dots) of the cone in Figure (ii).
\[
\epsfxsize = 10cm \epsfbox{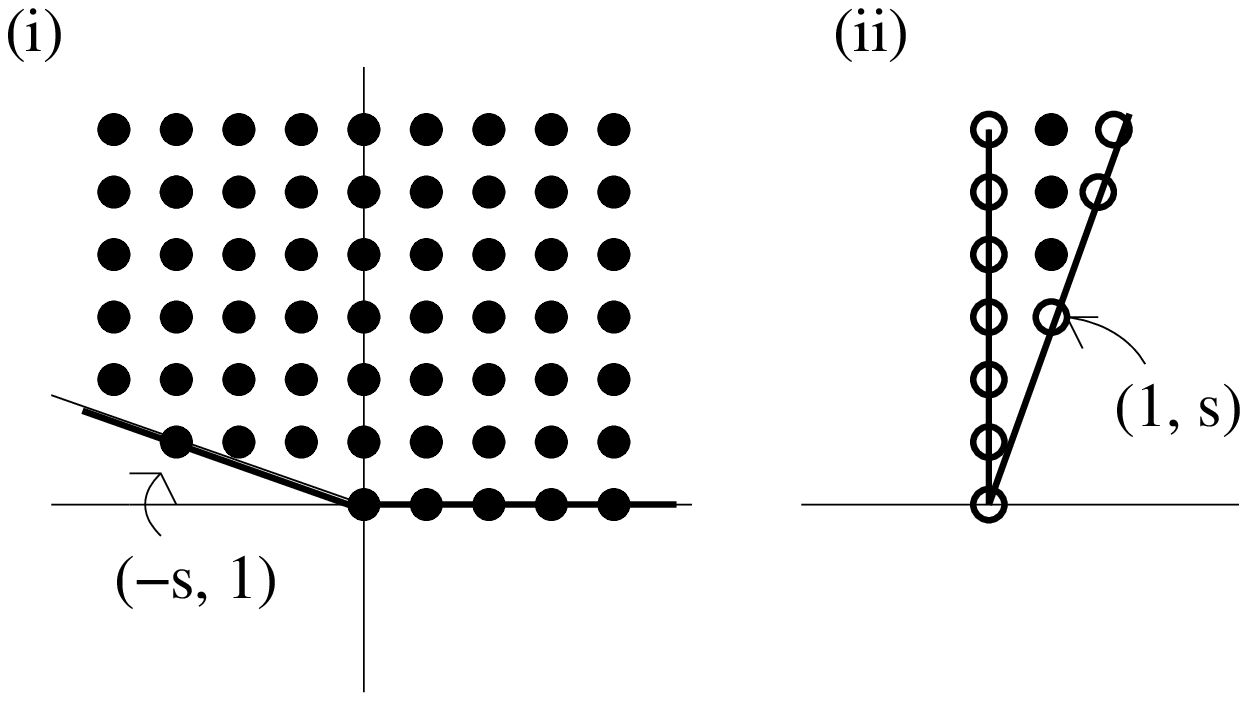}
\]

\begin{ex}\label{ex: bounds}
If we choose coarsening vector ${\bf v} = (1, s+1)$, then
$\vdeg(x_1)=\vdeg(x_3)=1$, $\vdeg(x_2)=1$ and $\vdeg(x_4)=s+1$. Therefore,
$c_{\bf v}=s+1$, $\vregnum(S)=2s$, and
$s_{\bf v}=3s$.

If we choose coarsening vector ${\bf u} = (1, 2s)$, then
$\udeg(x_1)=\udeg(x_3)=1$, $\udeg(x_2)=s$ and $\udeg(x_4)=2s$. Therefore,
 $c_{\bf u}=2s$, $\uregnum(S)=3s-1$, and
$s_{\bf u}=5s-2$.
In particular, when  $s=2$, then $\vregnum(S)=4$, and $\uregnum(S)=8$.
\end{ex}

\begin{ex}
Let $S$ be the coordinate ring of the Hirzebruch surface
$\mathbb{F}_s$, and consider the $S$-module $M = S/\langle
x_1x_2,x_3x_4\rangle$. This is the same ring as the ring of
Example \ref{ex3} except that we have renamed our indeterminates
and we have changed our grading.  So $H_{\mathfrak{m}}^i(M) = 0$
for all $i \neq 2$, and for $i = 2$, the fine $\NN^4$-graded
Hilbert series of $H_{\mathfrak{m}}^2(M)$ is the same as the one
given in Example \ref{ex3}. If we choose the positive coarsening
vector ${\bf v} = (1, s+1)$, then the $\ZZ$-graded Hilbert series
of $H_{\mathfrak{m}}^2(M^{[{\bf v}]})$ is
\[H(H_{\mathfrak{m}}^2(M^{[{\bf v}]});t) = 1 + 3(t^{-1}+t^{-2} + \cdots) + (t^{-s-1}+ t^{-2s-2} + \cdots) + \cdots;\]
in other words, $H_{\mathfrak{m}}^2(M)_p = 0$ for all $p \geq 1$ but $H_{\mathfrak{m}}^2(M)_0 \neq 0$.  So, $a^2(M^{[{\bf v}]}) = 0$
and $a^i(M^{[{\bf v}]}) = -\infty$ for $i \neq 2$.  Since $c_{\bf v} = (s+1)$ by Theorem
\ref{lem:a-inv} we have
\[\vregnum(M) = 0 - (s+1)(1-2) + 1 = s+2.\]
\end{ex}

\end{document}